# Complete Integral of Primer-Vector Equations for Transfers in a Central Gravitational Field

Sergey Zaborsky

RSC "Energia," Korolev, Moscow Region, 141070, Russia

Dr. Sergey Zaborsky

sergey.zaborsky@rsce.ru

**Abstract**

This paper demonstrates the existence of a complete integral for the system of differential equations of Lawden's primer-vector, which is used in the optimization of space transfers in a central gravitational field. The derived complete integral has been shown to significantly reduce the order of the differential system for the primer-vector from sixth to second, thereby simplifying the optimization problem into a boundary value problem with four parameters. The presence of a complete integral enables the exclusion of the transversality conditions, which introduce significant complexity to the boundary value problem. The problem of transfer optimization is considerably simplified due to the existence of the full integral and generating solutions. The analysis reveals that, depending on the given constraints, there are six types of optimization problems, each corresponding to a specific boundary value problem.

## 1 Introduction

In practical optimization of space maneuvers, nonlinear programming methods with equality and inequality constraints are commonly used [1-3]. These methods necessitate preliminary knowledge of maneuver parameters, such as the number of engine firings, their duration, and the intervals between them, critical for effective planning. Classical optimization methods, including the calculus of variations [4], dynamic programming [5], and methods based on Pontryagin's maximum principle [6], do not require such prior information. The calculus of variations demands continuity of the control function, which may limit its applicability in discontinuous scenarios, while dynamic programming can be impractical due to the need to process large volumes of information. Pontryagin's maximum principle offers a strategic advantage by reducing the complex problem of optimal control to a boundary value problem, importantly without requiring continuity in the control function. This principle involves integrating a system of six ordinary differential equations for variables that define the optimal control functions, such as Lawden's primer-vector [7].

This article demonstrates the existence of a complete integral for the primer-vector differential equation system, derived from a vector integral initially published by S. Pines [8]. The complete integral makes it possible to reduce the system of primer-vector differential equations from the sixth order to the second. "A long-standing challenge in optimal guidance for orbital insertion of launch vehicles and orbital transfers of spacecraft is how to reduce the transversality conditions in the optimal control problem." [9]. The presence of a complete integral enables the exclusion of the transversality conditions. Also, the existence of a complete integral when using generating solutions [10] significantly simplifies the problem of optimizing transfers. Depending on the specified constraints, six kinds of optimization problems for space maneuvers in a central gravitational field are identified, each corresponding to a unique boundary value problem.

## 2 Problem Statement and Derivation of the Complete Integral

The system of ordinary differential equations describing motion in a central gravitational field can be expressed in an arbitrary inertial Cartesian geocentric frame:

$$\begin{cases} \dfrac{d\Delta t}{dt} = 1 \\ \dfrac{d^2\boldsymbol{r}}{dt^2} + \dfrac{\gamma}{r^3}\boldsymbol{r} = \dfrac{d\Delta v}{dt}\boldsymbol{i} \\ \dfrac{d\Delta v}{dt} = \dfrac{P}{m_0}\exp\left(\dfrac{\Delta v}{g_0 I_{sp}}\right) \end{cases} \tag{1}$$

where $t$ represents time, $\boldsymbol{r}$ is the radius vector, $\gamma$ is the gravitational parameter, $P$ denotes engine thrust (which can either be zero or reach a maximum specified value $P_{\max}$), $\Delta t$ is the transfer time, $\Delta v$ is the change in characteristic velocity due to thrust, $\boldsymbol{i}$ is the unit vector in the direction of thrust acceleration, $I_{sp}$ is the propellant specific impulse (thrust/propellant weight flow weight) of the engine, $m_0$ represents the initial mass of the spacecraft, and $g_0$ is the gravitational constant at sea level.

The initial conditions at the start of the transfer specified at time $t_0$, as detailed in Eq. (2):

$$\boldsymbol{r} = \boldsymbol{r}_0 \qquad d\boldsymbol{r}/dt = \boldsymbol{v}_0 \qquad \Delta t_0 = 0 \qquad \Delta v_0 = 0 \tag{2}$$

The radius vector and velocity vector at the end of the transfer are specified at time $t_f$, as described in Eq. (3):

$$\boldsymbol{r} = \boldsymbol{r}_f \qquad d\boldsymbol{r}/dt = \boldsymbol{v}_f \tag{3}$$

The optimization problem is as follows: given the initial conditions Eq. (2), transform the Eqs. (1) using the control functions $P(t)$ and $\boldsymbol{i}(t)$. This transformation should ensure that the boundary conditions - represented by various functions of $\boldsymbol{r}_f$ and $v_f$ (e.g., target orbit parameters) - are satisfied in such a way that the condition of minimum total delta velocity $\Delta v_f$ from the thrust force is satisfied in Eq. (4):

$$\Delta v_{\Sigma}^{opt} = \min_{P(t),\boldsymbol{i}(t)} \Delta v_f \tag{4}$$

The control functions $P(t)$ and $\boldsymbol{i}(t)$ must be determined according to Pontryagin's maximum principle, which maximizes the Hamiltonian (5):

$$H = -\frac{\gamma}{r^3}\boldsymbol{r}\cdot\boldsymbol{p} - \boldsymbol{v}\cdot\frac{d\boldsymbol{p}}{dt} + \frac{P}{m_0}\exp\left(\frac{\Delta v}{g_0 I_{sp}}\right)\left(\psi_{\Delta v} + \boldsymbol{i}\cdot\boldsymbol{p}\right) + \psi_{\Delta t} \tag{5}$$

where $\boldsymbol{p}$ - primer-vector, $\boldsymbol{v} = d\boldsymbol{r}/dt$ - velocity vector, $\psi_{\Delta t}$ is conjugate variable associated with $\Delta t$ and $\psi_{\Delta v}$ is conjugate variable associated with $\Delta v$. The maximum of the Hamiltonian

$$\max_{P(t), \boldsymbol{i}(t)} H = -\frac{\gamma}{r^3}\boldsymbol{r}\cdot\boldsymbol{p} - \boldsymbol{v}\cdot\frac{d\boldsymbol{p}}{dt} + \frac{P}{m_0}\exp\left(\frac{\Delta v}{g_0 I_{sp}}\right)\kappa + \psi_{\Delta t} \tag{6}$$

is satisfied when

$$\boldsymbol{i} = \frac{\boldsymbol{p}}{p} \qquad P = \begin{cases} P_{\max}, \kappa \geq 0 \\ 0, \kappa < 0 \end{cases} \tag{7}$$

where $\kappa = \psi_{\Delta v} + p$ is the switching function. The primer-vector $\boldsymbol{p}$, which appears in the control functions $P(t)$ and $\boldsymbol{i}(t)$, is determined by solving the system of differential equations:

$$\begin{cases} \dfrac{d\psi_{\Delta t}}{dt} = -\dfrac{\partial H}{\partial \Delta t} = 0 \\ \dfrac{d^2\boldsymbol{p}}{dt^2} = \nabla H = \dfrac{3\gamma}{r^5}(\boldsymbol{r}\cdot\boldsymbol{p})\boldsymbol{r} - \dfrac{\gamma}{r^3}\boldsymbol{p} \\ \dfrac{d\psi_{\Delta v}}{dt} = -\dfrac{\partial H}{\partial \Delta v} = -\dfrac{1}{g_0 I_{sp}}\dfrac{d\Delta v}{dt}\kappa \end{cases} \tag{8}$$

The Hamiltonian (6) is equal to the first integral of the Eqs. (8) (herewith $\psi_{\Delta t} = const = C$ ), a concept noted in [7, p. 60, Eq. (3.39)], i.e. equal to 0:

$$-\frac{\gamma}{r^3}\boldsymbol{r}\cdot\boldsymbol{p} - \boldsymbol{v}\cdot\frac{d\boldsymbol{p}}{dt} + \frac{P}{m_0}\exp\left(\frac{\Delta v}{g_0 I_{sp}}\right)\kappa + C = 0 \tag{9}$$

By forming the vector cross product of $\boldsymbol{p}$ with second Eq. (1), the vector cross product of $\boldsymbol{r}$ with second Eq. (8), and adding them, we obtain the following equation $\boldsymbol{p}\times d^2\boldsymbol{r}/dt^2 + \boldsymbol{r}\times d^2\boldsymbol{p}/dt^2 = \boldsymbol{0}$ having vector integral [8]:

$$\frac{d\boldsymbol{p}}{dt}\times\boldsymbol{r} - \boldsymbol{p}\times\boldsymbol{v} - \boldsymbol{Z} = \boldsymbol{0} \tag{10}$$

where $\boldsymbol{Z} = \boldsymbol{const}$ is a constant vector.

By forming the dot product of Eq. (10) with the radius unit vector $\boldsymbol{r}$, with the transverse unit vector $\boldsymbol{h}\times\boldsymbol{r}$ and with the angular momentum unit vector $\boldsymbol{h} = \boldsymbol{r}\times\boldsymbol{v}$, we obtain with vector algebra [11, Chap. 5.2] and having in mind that $r^2/h = dt/d\theta$ ($\theta$ is the angular flight distance)

$$\begin{bmatrix} \frac{\boldsymbol{r}}{r} \\ \frac{\boldsymbol{h}\times\boldsymbol{r}}{hr} \\ \frac{\boldsymbol{h}}{h} \end{bmatrix} \cdot \left( \frac{d\boldsymbol{p}}{dt}\times\boldsymbol{r} - \boldsymbol{p}\times\boldsymbol{v} - \boldsymbol{Z} \right) = \begin{bmatrix} \frac{1}{r}\left(\boldsymbol{h}\cdot\boldsymbol{p} - \boldsymbol{r}\cdot\boldsymbol{Z}\right) \\ \frac{1}{hr}\left( \left(\boldsymbol{h}\cdot\frac{d\boldsymbol{p}}{dt}\right) r^2 - \left(\boldsymbol{h}\cdot\boldsymbol{r}\right)\left(\boldsymbol{r}\cdot\frac{d\boldsymbol{p}}{dt}\right) - \left(\boldsymbol{r}\cdot\boldsymbol{v}\right)\left(\boldsymbol{h}\cdot\boldsymbol{p}\right) + \left(\boldsymbol{h}\cdot\boldsymbol{v}\right)\left(\boldsymbol{r}\cdot\boldsymbol{p}\right) - \left(r^2\boldsymbol{v} - \left(\boldsymbol{r}\cdot\boldsymbol{v}\right)\boldsymbol{r}\right)\cdot\boldsymbol{Z} \right) = \\ \frac{r^2}{h}\left( \frac{1}{r}\frac{d\left(\boldsymbol{h}\cdot\boldsymbol{p}\right)}{dt} - \frac{\boldsymbol{r}\cdot\boldsymbol{v}}{r^3}\left(\boldsymbol{h}\cdot\boldsymbol{p}\right) - \left(\frac{\boldsymbol{v}}{r} - \frac{\boldsymbol{r}\cdot\boldsymbol{v}}{r^3}\boldsymbol{r}\right)\cdot\boldsymbol{Z} \right) = \frac{dt}{d\theta}\left( \frac{d}{dt}\left(\frac{\boldsymbol{h}\cdot\boldsymbol{p}}{r}\right) - \frac{d}{dt}\left(\frac{\boldsymbol{r}\cdot\boldsymbol{Z}}{r}\right) \right) = \\ \frac{d}{d\theta}\left( \frac{1}{r}\left(\boldsymbol{h}\cdot\boldsymbol{p} - \boldsymbol{r}\cdot\boldsymbol{Z}\right) \right) \\ \frac{1}{h}\left( \left(\boldsymbol{h}\times\boldsymbol{v}\right)\cdot\boldsymbol{p} - \left(\boldsymbol{h}\times\boldsymbol{r}\right)\cdot\frac{d\boldsymbol{p}}{dt} - \boldsymbol{h}\cdot\boldsymbol{Z} \right) \end{bmatrix} = \boldsymbol{0} \tag{11}$$

The first expression in Eq. (11) is a complete integral of the second equation for the primer-vector with respect to the angular flight distance. Thus, the second expression in Eq. (11) can be disregarded, facilitating a reduction in the order of the differential system for the primer-vector from sixth to second order, as demonstrated in:

$$\begin{cases} -\frac{\gamma}{r^3}\boldsymbol{r}\cdot\boldsymbol{p} - \boldsymbol{v}\cdot\frac{d\boldsymbol{p}}{dt} + \frac{P}{m_0}\exp\left(\frac{\Delta v}{g_0 I_{sp}}\right)\kappa + C = 0 \\ \boldsymbol{h}\cdot\boldsymbol{p} = \boldsymbol{r}\cdot\boldsymbol{Z} \\ \left(\boldsymbol{h}\times\boldsymbol{v}\right)\cdot\boldsymbol{p} - \left(\boldsymbol{h}\times\boldsymbol{r}\right)\cdot\frac{d\boldsymbol{p}}{dt} = \boldsymbol{h}\cdot\boldsymbol{Z} \end{cases} \tag{12}$$

Let us assume that the projection of vector $\boldsymbol{Z}$ onto the direction $\boldsymbol{r}_0$ is equal to a constant $E$, and its projection onto the direction $\boldsymbol{h}_0$ is equal to a constant $A$. Given that the third condition can be ignored, the projection of vector $\boldsymbol{Z}$ onto the transverse direction $\boldsymbol{h}_0 \times \boldsymbol{r}_0$ can be set arbitrarily, including to zero. This simplification leads to $\boldsymbol{Z}$ depending only on two constants $E$ and $A$, as outlined in Eq. (13):

$$\boldsymbol{Z} = E\frac{\boldsymbol{r}_0}{r_0} + A\frac{\boldsymbol{h}_0}{h_0} \tag{13}$$

After substituting Eq. (13) into Eq. (12), the system of sixth-order equations for the primer-vector is replaced with the simplified form:

$$
\begin{cases}
-\dfrac{\gamma}{r^3}\boldsymbol{r}\cdot\boldsymbol{p}-\boldsymbol{v}\cdot\dfrac{d\boldsymbol{p}}{dt}+\dfrac{d\Delta v}{dt}\kappa+C=0\\
\boldsymbol{h}\cdot\boldsymbol{p}=Er\cos\theta\\
(\boldsymbol{h}\times\boldsymbol{v})\cdot\boldsymbol{p}-(\boldsymbol{h}\times\boldsymbol{r})\cdot\dfrac{d\boldsymbol{p}}{dt}=Ah\cos\chi\\
\dfrac{d\psi_{\Delta v}}{dt}=-\dfrac{1}{g_0 I_{sp}}\dfrac{d\Delta v}{dt}\kappa
\end{cases}
\tag{14}
$$

Where $\cos\theta=\boldsymbol{r}_0\cdot\boldsymbol{r}/(r_0 r)$, and $\chi$ represents the angle between the planes of the initial and current orbits ( $\cos\chi=\boldsymbol{h}_0\cdot\boldsymbol{h}/(h_0 h)$ ).

## 3 Problem Statement in an Orbital Non-Inertial Frame

To describe motion in an orbital non-inertial frame (labelled with the symbol "^"), it is necessary to make the following substitution in the Eqs. (1)

$$
\boldsymbol{r}\to\hat{\boldsymbol{r}}=\begin{bmatrix} r\\ 0\\ \hat{z}=0\end{bmatrix}\qquad \frac{d^2\boldsymbol{r}}{dt^2}\to\frac{d^2\hat{\boldsymbol{r}}}{dt^2}+\frac{\hat{\boldsymbol{h}}}{r^2}\times\hat{\boldsymbol{v}}\qquad \boldsymbol{i}\to\hat{\boldsymbol{i}}=\frac{\hat{\boldsymbol{p}}}{\hat{p}}=\begin{bmatrix} S^0\\ T^0\\ W^0\end{bmatrix}
\tag{15}
$$

and in the Eqs. (8) [7, p.82, Eq. (5.16)]

$$
\boldsymbol{p}\to\hat{\boldsymbol{p}}=\begin{bmatrix}\lambda\\ \mu\\ \nu\end{bmatrix}\qquad \frac{d\boldsymbol{p}}{dt}\to\frac{d\hat{\boldsymbol{p}}}{dt}+\frac{\hat{\boldsymbol{h}}}{r^2}\times\hat{\boldsymbol{p}}\qquad \kappa\to\hat{\kappa}=\psi_{\Delta v}+\hat{p}
\tag{16}
$$

Considering that

$$
\hat{\boldsymbol{v}}=\begin{bmatrix} v_r\\ v_\theta\\ v_{\hat{z}}=0\end{bmatrix}\qquad \hat{\boldsymbol{h}}=\begin{bmatrix}0\\ 0\\ rv_\theta\end{bmatrix}\qquad \frac{d^2\hat{z}}{dt^2}=\frac{dv_{\hat{z}}}{dt}=\frac{d\Delta v}{dt}W^0
\tag{17}
$$

where $v_r$, $v_\theta$ and $v_{\hat{z}}$ are the radial, transversal, and normal components of the velocity vector $\hat{\boldsymbol{v}}$, accordingly, and $S^0$, $T^0$ and $W^0$ are the radial, transversal, and normal of the unit vector $\hat{\boldsymbol{i}}$ , accordingly, the system of equations in the orbital non-inertial frame is written in the following expanded form (in this case from Eq. $dv_{\hat{z}}=v_\theta d\chi$ follows Eq. $d\chi/dt=d\Delta v/dt\,W^0/v_\theta$ .)

$$
\begin{cases}
\dfrac{d\Delta t}{dt} = 1 \\
\dfrac{dr}{dt} = v_r \\
\dfrac{dv_r}{dt} = \left(v_\theta^2 - \dfrac{\gamma}{r}\right)\dfrac{1}{r} + \dfrac{d\Delta v}{dt} S^0 \\
\dfrac{dv_\theta}{dt} = -\dfrac{v_r v_\theta}{r} + \dfrac{d\Delta v}{dt} T^0 \\
\dfrac{d\chi}{dt} = \dfrac{1}{v_\theta}\dfrac{d\Delta v}{dt} W^0 \\
\dfrac{d\theta}{dt} = \dfrac{v_\theta}{r} \\
\dfrac{d\Delta v}{dt} = \dfrac{P}{m_0}\exp\left(\dfrac{\Delta v}{g_0 I_{sp}}\right)
\end{cases}
\tag{18}
$$

and Eqs. (14) in the form

$$
\begin{cases}
-\dfrac{\gamma}{r^2}\lambda - v_r\left(\dfrac{d\lambda}{dt} - \dfrac{v_\theta}{r}\mu\right) - v_\theta\left(\dfrac{d\mu}{dt} + \dfrac{v_\theta}{r}\lambda\right) + \dfrac{d\Delta v}{dt}\hat{\kappa} + C = 0 \\
v_\theta \nu = E \Rightarrow \dfrac{d\nu}{dt} = \dfrac{v_r}{r}\nu \Rightarrow \dfrac{d^2\nu}{dt^2} = \left(-\dfrac{\gamma}{r^3} + \left(\dfrac{v_\theta}{r}\right)^2\right)\nu \\
v_r\mu - v_\theta\lambda - r\left(\dfrac{d\mu}{dt} + \dfrac{v_\theta}{r}\lambda\right) = A \\
\dfrac{d\psi_{\Delta v}}{dt} = -\dfrac{1}{g_0 I_{sp}}\dfrac{d\Delta v}{dt}\hat{\kappa}
\end{cases}
\tag{19}
$$

Solving Eqs. (19) with respect to $d\lambda/dt$ and $d\mu/dt$, Eq. (8) can be written as follows

$$
\begin{cases}
\dfrac{d\lambda}{dt} = \dfrac{1}{v_r}\left(\left(v_\theta^2 - \dfrac{\gamma}{r}\right)\dfrac{1}{r}\lambda + \dfrac{v_\theta}{r}A + C + \dfrac{d\Delta v}{dt}\hat{\kappa}\right) \\
\dfrac{d\mu}{dt} = \dfrac{1}{r}\left(v_r\mu - 2v_\theta\lambda - A\right) \\
\nu = \dfrac{1}{v_\theta}E \\
\dfrac{d\psi_{\Delta v}}{dt} = -\dfrac{1}{g_0 I_{sp}}\dfrac{d\Delta v}{dt}\hat{\kappa}
\end{cases}
\tag{20}
$$

Initial conditions for integrating Eqs. (18) and (20) at time $t_0$

$$
r = r_0 \qquad v_r = v_{0r} \qquad v_\theta = v_{0\theta} \qquad \theta_0 = 0 \quad \chi_0 = 0 \quad \Delta v_0 = 0 \qquad \psi_{\Delta v0} = -1 \tag{21}
$$

Analysis of the Eqs. (20) reveals that if $E = 0$, the primer-vector lies within the orbital plane, making the transfer coplanar ($\chi = 0$ and $\Delta\chi = \cos^{-1}\left(\boldsymbol{h}_0 \cdot \boldsymbol{h}_f / \left(h_0 h_f\right)\right) = 0$). It is shown in [7] that for $C \neq 0$, the flight time $\Delta t$ must be specified. In [12], in the analysis of Lawden's solutions [7], it is shown that when $C = 0$ and $A = 0$ both the flight time $\Delta t = t_f - t_0$ and the transfer angle $\Delta\theta = \cos^{-1}\left(\boldsymbol{r}_0 \cdot \boldsymbol{r}_f / \left(r_0 r_f\right)\right)$ are unspecified and optimized. Accordingly, if $\Delta\theta$ is set and $\Delta t$ is not set, then $A \neq 0$ and $C = 0$. The boundary conditions and the corresponding parameters, which are purposefully varied when solving the boundary value problem, are shown in Table 1.

**Table 1 Classification of minimization problems $\Delta v_f$**

| Kind | | Variable parameters | | | | Boundary conditions | | | | |
|---|---|---|---|---|---|---|---|---|---|---|
| I | Coplanar transfer $(E = 0)$ | $\lambda_0$ | $\mu_0$ | - | - $(A = C = 0)$ | $v_{fr}$ | $v_{f\theta}$ | - | - | $r_f$ |
| | Noncoplanar transfer $(E \neq 0)$ | $\lambda_0$ | $\mu_0$ | $E$ | - $(A = C = 0)$ | $v_{fr}$ | $v_{f\theta}$ | $\chi_f$ | - | |
| II | Coplanar transfer $(E = 0)$ | $\lambda_0$ | $\mu_0$ | - | $A$ $(C = 0)$ | $v_{fr}$ | $v_{f\theta}$ | - | $r_f$ | $\Delta\theta$ |
| | Noncoplanar transfer $(E \neq 0)$ | $\lambda_0$ | $\mu_0$ | $E$ | $A$ $(C = 0)$ | $v_{fr}$ | $v_{f\theta}$ | $\chi_f$ | $r_f$ | |
| III | Coplanar transfer $(E = 0)$ | $\lambda_0$ | $\mu_0$ | - | $C$ $(A = 0)$ | $v_{fr}$ | $v_{f\theta}$ | - | $r_f$ | $\Delta t$ |
| | Noncoplanar transfer $(E \neq 0)$ | $\lambda_0$ | $\mu_0$ | $E$ | $C$ $(A = 0)$ | $v_{fr}$ | $v_{f\theta}$ | $\chi_f$ | $r_f$ | |

The choice of one of the boundary conditions as the end condition for the integration of differential Eqs. (18) leads to the equality of the number of variable parameters to the number of boundary conditions. For example, for kind I it may be $r_f$, for kind II - $\Delta\theta$ and for kind III - $\Delta t$ (last column in Table 1). The optimization problem is reduced to the following boundary value problem: to find such initial elements of the primer-vector $\lambda_0$ and $\mu_0$, constants $E$ and $A$ or $C$, so that as a result of the joint integration of the Eqs. (18) and Eqs. (20), any kind of boundary conditions from Table 1 is satisfied. If the elements of the orbit are chosen as boundary conditions, then one can use the system of equations of the rates of the six classical elements [13, p. 200, Eq. (9.21)] instead of the system of Eqs. (18).

When $\kappa = 0$, the analytical solutions of the equations for $\lambda$ and $\mu$ in (20) coincide with the Louden solutions [7, p. 85, (5.55), (5.56)]. It should be noted that the Lawden's solution for $\nu$ [7, p. 86, Eq. (5.57)] used in primer-vector researches, e.g., [14] and [15], does not satisfy the Pines's vector integral (10), the complete integral (20), and the solution of the impulsive noncoplanar transfers optimization problem [16, Eq. (11)]. This is because Pines's vector integral was published later than Lawden's results and could not be considered by him. Lawden did not take into account that in [7, p. 82, Eqs. (5.24) and (5.25)] there may be a case when $v/r = const$, following from Eq. (20).

## 4 Numerical Examples

The optimal transfer between two given orbits with two engine burns is considered. The results of solving for constraints of kinds I, II and III are presented. The integration of the Eqs. (18) and Eqs. (20) was performed using the 4th-order Runge-Kutta method with a step size of 1.0 s. Initial approximations for the unknown parameters in boundary value optimization problems can be obtained using the generating solutions method [10]. It is necessary to determine the time points when the switching function $\kappa = 0$ with the greatest possible accuracy. In the process of solving the boundary value problem there may occur cases when $\kappa < 0$ at the initial moment of movement. In this case, it is recommended to move the initial point of the movement "backwards". The behavior of the function $\kappa$ at various values is illustrated in numerical examples [10]. It can be concluded that the engine burn at the beginning of the transfer is not always optimal. During the solution process, various schemes for switching on the propulsion system may appear. It is also necessary to analyze the resulting local extremes.

Table 2 provides the initial data used in the numerical example.

**Table 2 Constants**

| $P_{\max}$, N | $I_{sp}$, s | $m_0$, кg | $\gamma$, m$^3$/s$^2$ | $g_0$, m/s$^2$ |
|---|---|---|---|---|
| 74270 $g_0$ | 356 | $28 \cdot 10^3$ | $398600.436233 \cdot 10^9$ | 9.80665 |

Table 3 presents the initial conditions of the motion and parameters of the initial orbit ( $l$ - semi-latus rectum, $e$ - eccentricity, $i$ - inclination, $f$ - true anomaly).

**Table 3 Initial conditions**

| $r_0$, km | $v_{0r}$, m/s | $v_{0\theta}$, m/s | $l_0$, km | $e_0$ | $i_0$, deg | $f_0$, deg |
|---|---|---|---|---|---|---|
| 6553.71 | 74.57 | 6994.07 | 5271.04 | 0.195904 | 0.0 | 177.49 |

Table 4 provides the boundary conditions and parameters of the final orbit used in the numerical example.

**Table 4 Final conditions**

| $r_f$, km | $v_{fr}$, m/s | $v_{f\theta}$, m/s | $l_f$, km | $e_f$ | $i_f$, deg | $f_f$, deg |
|---|---|---|---|---|---|---|
| 11595.00 | 2913.68 | 6685.04 | 15073.42 | 0.641120 | $0 \div 90$ | 62.101 |

The generating solutions of solving for constraints of kind I ( $A = C = 0$ ) and $P_{\max} \to \infty$ ( $\hat{p}_0 = 1$, $\hat{\kappa}_0 = 0$ and $E = \partial \Delta v_{\Sigma}^{opt} / \partial \Delta \chi$ ) are presented in Table 5.

**Table 5 Generating Solutions of the boundary value problem** ( $P_{\max} \to \infty$ )

| $i_f = \Delta\chi$, deg | $\Delta v_1$, m/s | $\Delta v_2$, m/s | $\Delta v_{\Sigma}^{opt}$, m/s | $\lambda_0$ | $\mu_0$ | $E$, m/s | $\theta_f$, deg | $t_f$, s |
|---|---|---|---|---|---|---|---|---|
| 0 | 2389.70 | 1405.38 | 3795.08 | 0.006226 | 0.999980 | 0.0 | 112.16 | 2173.62 |
| 10 | 2465.52 | 1523.86 | 3989.38 | 0.006692 | 0.973021 | 2160.83 | 112.71 | 2189.45 |
| 20 | 2627.49 | 1884.57 | 4512.06 | 0.008047 | 0.917763 | 3703.34 | 114.45 | 2239.85 |
| 60 | 2607.57 | 5192.90 | 7800.47 | 0.018077 | 0.839295 | 4910.80 | 131.55 | 2739.63 |
| 90 | 2283.39 | 7928.60 | 10211.99 | 0.025574 | 0.881815 | 4202.54 | 143.30 | 3100.35 |

The results of solving for constraints of kind I ( $A = C = 0$ ) are presented in Table 6.

**Table 6 Solutions of the boundary value problem of kind I**

| $i_f = \Delta\chi$, deg | $\Delta v_1$, m/s | $\Delta v_2$, m/s | $\Delta v_{\Sigma}^{opt}$, m/s | $\lambda_0$ | $\mu_0$ | $E$, m/s | $\theta_{f1}$, deg | $t_{f1}$, s |
|---|---|---|---|---|---|---|---|---|
| 0 | 2384.70 | 1411.24 | 3795.94 | 0.003123 | 0.9999777 | 0.0 | 116.70 | 2242.29 |
| 10 | 2455.17 | 1535.31 | 3990.48 | 0.003351 | 0.9534582 | 2109.13 | 117.06 | 2256.63 |
| 20 | 2578.87 | 1937.07 | 4515.94 | 0.003014 | 0.871514 | 3429.72 | 119.46 | 2321.09 |
| 60 | 2540.29 | 5409.55 | 7949.84 | 0.001956 | 0.770359 | 4459.84 | 135.70 | 2778.10 |
| 90 | 2451.06 | 8380.07 | 10831.13 | 0.020320 | 0.761445 | 4534.38 | 145.03 | 3053.78 |

The results of solving for constraints of kind II ($C = 0$, $\Delta\theta = \theta_{fl} - 5\deg$) are presented in Table 7.

**Table 7 Solutions of the boundary value problem of kind II**

| $i_f = \Delta\chi$, deg | $\Delta v_1$, m/s | $\Delta v_2$, m/s | $\Delta v_\Sigma^{opt}$, m/s | $\lambda_0$ | $\mu_0$ | $E$, m/s | $A$ m/s | $\Delta\theta$, deg |
|---|---|---|---|---|---|---|---|---|
| 0 | 2419.42 | 1383.65 | 3803.07 | 0.078068 | 1.0096411 | 0.0 | -184.57 | 111.70 |
| 10 | 2496.55 | 1501.59 | 3998.14 | 0.007974 | 0.9623180 | 2136.57 | -184.95 | 112.06 |
| 20 | 2627.50 | 1893.90 | 4521.40 | 0.068049 | 0.8773187 | 3473.83 | -146.85 | 114.46 |
| 60 | 2586.43 | 5368.45 | 7954.88 | 0.041325 | 0.7725987 | 4491.76 | -67.80 | 130.70 |
| 90 | 2487.01 | 8351.57 | 10838.58 | 0.035000 | 0.7636846 | 4556.92 | -51.91 | 140.03 |

The results of solving for constraints of kind III ($A = 0$, $\Delta t = t_{fl} - 35s$) are presented in Table 8.

**Table 8 Solutions of the boundary value problem of kind III**

| $i_f = \Delta\chi$, deg | $\Delta v_1$, m/s | $\Delta v_2$, m/s | $\Delta v_\Sigma^{opt}$, m/s | $\lambda_0$ | $\mu_0$ | $E$, m/s | $C$, m/s$^2$ | $\Delta t$, deg |
|---|---|---|---|---|---|---|---|---|
| 0 | 2401.78 | 1395.02 | 3796.80 | 0.020752 | 1.0032990 | 0.0 | -0.05988 | 2277.29 |
| 10 | 2477.44 | 1514.02 | 3991.46 | 0.022541 | 0.9565211 | 2136.50 | -0.06476 | 2221.63 |
| 20 | 2603.84 | 1912.41 | 4516.25 | 0.018168 | 0.8728545 | 3451.67 | -0.04750 | 2286.09 |
| 60 | 2558.14 | 5392.82 | 7950.96 | 0.008729 | 0.7703796 | 4471.14 | -0.01589 | 2743.10 |
| 90 | 2463.51 | 8369.72 | 10833.23 | 0.007036 | 0.7615636 | 4541.35 | -0.01066 | 3018.78 |

Analyzing data from the Tables 6-8, we can conclude that the most optimal kind is I.

## 5 Conclusion

This paper establishes the existence of a complete integral for the primer-vector system in optimizing transfers in a central gravitational field. By reducing the complexity of the differential system from sixth to second order, the paper clarifies how such reductions simplify the optimization process to a four-parameter boundary value problem. It should be noted that the obtained results allow us to solve the problem of minimizing the flight time. In this case, if, when solving the boundary value problem of kind III, we set the boundary value $\Delta v$ instead of $\Delta t$, then instead of $\Delta v_f$, $\Delta t_f$ will be minimized.